\documentclass[11pt]{article}
\usepackage{latexsym}
\usepackage{amscd}
\usepackage{amssymb}
\usepackage{amsmath}
\usepackage{amsthm}
\usepackage{enumerate}
\usepackage{verbatim}
\def\co{\colon\thinspace}

\newtheorem{thm}{Theorem}[section]

\newtheorem{lem}[thm]{Lemma}

\newtheorem{Example}[thm]{Example}

\newtheorem{remark}[thm]{Remark}
\newenvironment{rmk}{\begin{remark}\rm}{\end{remark}}
\newtheorem{Fact}[thm]{Fact}

\newtheorem{Nothing}[thm]{$\!\!\!$}

\begin{document}
\abovedisplayskip=6pt plus3pt minus3pt
\belowdisplayskip=6pt plus3pt minus3pt
\title{\bf 
Pinching surface groups in \\
complex hyperbolic plane}
\author{Igor Belegradek
\ \\
\\ fax: $(626)$ $585-1728$\\
{\it email: ibeleg@its.caltech.edu}\\
Department of Mathematics 253-37\\
California Institute of Technology\\
Pasadena, CA, $91125$, USA\\
\ \\
{\it keywords:} complex hyperbolic geometry, surface group.
\ \\
{\it Dedicated to the memory of Hanna Sandler}}
\date{}
\maketitle
\begin{abstract} 
We construct first examples of discrete
geometrically finite subgroups of $\mathrm{PU}(2,1)$
which contain parabolic elements, and are isomorphic
to surface groups of genus $\ge 2$. 
\end{abstract}

\section{Introduction}
Even though complex hyperbolic Kleinian groups 
lack the flexibility
of their real hyperbolic cousins,
they do come in many shapes and sizes,
and there is no structure theory in sight.
At this stage it seems useful to work out individual 
examples, such as the ones considered in this paper.
 
Let $M$ be a closed orientable surface of genus $g\ge 2$.
In this note we construct new examples of
discrete embeddings of $\pi=\pi_1(M)$
into $\mathrm {PU}(2,1)$, the full group of biholomorphic isometries
of the complex hyperbolic plane $\mathbf {H}^2_{\mathbb  C}$.
The group $\pi$ can be realized as a lattice
in the subgroups $\mathrm {SO}(2,1)$ and $\mathrm {SU}(1,1)$ of 
$\mathrm {PU}(2,1)$, so there are two obvious discrete embeddings 
$\rho_r$ and $\rho_c$ of $\pi$ into $\mathrm{PU}(2,1)$.
The group $\rho_r(\pi)$ stabilizes a totally real plane 
in $\mathbf {H}^2_{\mathbb C}$, and the quotient complex hyperbolic 
surface $\mathbf {H}^2_{\mathbb C}/\rho_r(\pi)$ 
is the total space of the 
tangent bundle to $M$.
Similarly, $\rho_c(\pi)$ preserves a complex geodesic, and 
$\mathbf {H}^2_{\mathbb C}/\rho_c(\pi)$ 
is the square root of the tangent bundle to $M$~\cite{GKL}.

These two representations are also distinguished by 
the so-called Toledo invariant $\tau$,
which associates to a representation
$\rho\in\mathrm{Hom}(\pi_1(M),\mathrm{PU}(2,1))$ 
the (normalized) integral over $M$ of the pullback of the
K\"ahler form 
via a section of the flat $\mathbf {H}^2_{\mathbb C}$-bundle
corresponding to $\rho$.
In fact, $\tau(\rho_r)=0$ and $\tau(\rho_c)=\pm\chi(M)$.
According to~\cite{Tol, GKL}, $\tau$ 
is a $\frac{2}{3}\mathbb Z$-valued
locally constant function on the representation
space $\mathrm{Hom}(\pi_1(M),\mathrm{PU}(2,1))$ 
satisfying $|\tau|\le 2g-2$. E.~Xia
showed~\cite{Xia} that $\tau$ classifies the connected
components of $\mathrm{Hom}(\pi_1(M),\mathrm{PU}(2,1))$, and
D.~Toledo proved that $|\tau(\rho)|=2g-2$ iff
$\rho$ is an isomorphism onto a cocompact lattice
in the stabilizer of a complex geodesic.

By amalgamating the representations $\rho_r,\rho_c$,
W.~Goldman, M.~Kapovich, and B.~Leeb~\cite{GKL} showed that
each even value of $\tau$ is realized by a  
faithful discrete representations $\rho$ such that the complex
hyperbolic surface $\mathbf {H}^2_{\mathbb C}/\rho(\pi)$
is an oriented $\mathbb R^2$-bundle
over $M$ with Euler number $2g-2+|\tau(\rho)/2|$.

In all the above examples, the group $\rho(\pi)$ is 
geometrically finite without parabolics.
Goldman asked if this is always the case
for faithful discrete representations.
We provide a negative answer as follows. 

\begin{thm} \label{main thm}
Let $M$ be a closed orientable surface
of genus $g\ge 2$, and let $\gamma\in\pi_1(M)=\pi$ be a nontrivial 
element represented by a simple closed curve that separates $M$.
Then there exists a faithful discrete representation
$\rho\co\pi\to\mathrm {PU}(2,1)$ such that
$\rho(\pi)$ is geometrically finite, any maximal parabolic subgroup
of $\rho(\pi)$ is generated by a conjugate of $\rho(\gamma)$, the
quotient $\mathbf {H}^2_{\mathbb C}/\rho(\pi)$
is diffeomorphic to the tangent bundle of $M$, and $\tau(\rho)=0$.
\end{thm}

Loosely speaking, any nontrivial element $\gamma\in\pi_1(M)$ 
represented by a simple closed curve that separates $M$ can be 
pinched (i.e. made parabolic)
by some faithful discrete representation $\rho$. 
The element $\rho(\gamma)$ is 
called an {\it accidental parabolic}. 
A simple modification of our construction
yields faithful discrete representations
with several conjugacy classes of 
accidental parabolics, however the result we get is not
optimal, so do not write down the details.

We do not know which nonzero 
values of the Toledo invariant can be realized
by faithful discrete representations with accidental parabolics.
As we mentioned above, the components of 
$\mathrm{Hom}(\pi_1(M),\mathrm{PU}(2,1))$
with $|\tau|=2g-2$ cannot contain such representations.

The group $\rho(\pi)$ in Theorem~\ref{main thm} 
is obtained by amalgamating two discrete groups, each being 
a noncocompact lattice
in the stabilizer of some totally real plane,
along a common cyclic parabolic subgroup. 
Discreteness of the amalgamated product
is proved using the topological version of 
the Maskit combination theorem.

In author's opinion this method of proving
discreteness has some advantages over
the complex hyperbolic version of the Poincar\'e 
polyhedron theorem which was recently developed in~\cite{FZ, GP}.
Namely, the Poincar\'e theorem usually requires
explicit knowledge of the geometry of the 
fundamental polyhedron. 
By contrast, the Maskit combination 
theorem is stated in purely topological terms, in particular,
one does not need the fundamental polyhedron
to be bounded by bisectors, or any other
special hypersurfaces. This soft nature
of the combination theorems makes them easier to use.

This paper is a revised version of the preprint~\cite{bel}
written in 1995 when the author was a student at 
the University of Maryland.
It is a pleasure to thank Bill Goldman for numerous
helpful discussions, and Robert Miner for comments on the
first version of the paper.

\section{Vertical and horizontal translations in
the Siegel domain}
\label{sec: vh}
In this section we set up notations and collect some elementary
facts about geometry of the complex hyperbolic plane. 
The reader is referred to~\cite{Gol} for more information.
One of the standard models for the complex hyperbolic plane is 
the Siegel domain 
\[\mathbf {H}^2_{\mathbb C}=\{(w_1,w_2)\in 
\mathbb C^2\co  w_1\overline{w_1}<w_2+\overline{w_2}\}.\]
The real hypersurface 
\[
\{(w_1,w_2)\in \mathbb C^2\co  w_1\overline{w_1}=w_2+\overline{w_2}\}
\]
corresponds to the sphere at infinity with one point removed; 
the point is denoted $\infty $.
We identify the stabilizer of $\infty$ in 
$\mathrm{PU}(2,1)$ with the Heisenberg group $\mathcal H$.
One can introduce horospherical coordinates in $\mathbf {H}^2_{\mathbb C}$
by \[(x, y, u, v)=(z, u, v)\in 
\mathbb C\times (0,\infty)\times \mathbb R\]
where $z=x+iy=w_1$ and $u+iv=2\overline{w_2}-w_1\overline{w_1}$.

For every fixed $u>0$, the real hypersurface
$\{(z, u, v)\co z\in\mathbb C, v\in\mathbb R\}$ is a horosphere centered at 
the point $\infty$.
The group $\mathcal H$ acts simply
transitively on each horosphere, so that
every horosphere gets a Heisenberg space structure given by
\[(z_1, u, v_1)+(z_2, u, v_2)=
(z_1+z_2, u, v_1+v_2+2\mathrm {Im}(z_1\overline {z_2})).\]
Similarly, the hypersurface $u=0$, that corresponds to 
$\partial_\infty\mathbf {H}^2_{\mathbb C}\setminus\{\infty\}$,
has a simply transitive $\mathcal H$-action, and 
the $\mathcal H$-action on $Y=\mathbf {H}^2_{\mathbb C}\cup
\partial_\infty\mathbf {H}^2_{\mathbb C}\setminus\{\infty\}$
is smooth, free, and proper.

The real hyperbolic plane $\mathbf  {H}^2_{\mathbb R}$ sits inside 
$\mathbf  {H}^2_{\mathbb C}$ as a totally real $2$-plane
$\{(x, 0, u, 0)\in\mathbf  {H}^2_{\mathbb C}\}$.
(In this paper we always think of $\mathbf  {H}^2_{\mathbb R}$ as an 
upper half plane and think 
of $\mathbf  {H}^2_{\mathbb C}$ as the Siegel domain).
The orthogonal projection $\mathbf  {H}^2_{\mathbb C}\to 
\mathbf  {H}^2_{\mathbb R}$ is a $\mathbf  {SO}(2,1)$-equivariant smooth
$2$-plane bundle where 
the fibers are totally real $2$-planes
orthogonal to $\mathbf  {H}^2_{\mathbb R}$.
The projection extends to a $\mathbf  {SO}(2,1)$-equivariant smooth
map $\mathbf {H}^2_{\mathbb C}\cup
\partial_\infty\mathbf {H}^2_{\mathbb C}\to \mathbf {H}^2_{\mathbb R}\cup
\partial_\infty\mathbf {H}^2_{\mathbb R}$ which is the identity
on $\partial_\infty\mathbf {H}^2_{\mathbb R}$ and is a
closed disk bundle away from $\partial_\infty\mathbf {H}^2_{\mathbb R}$.
Since any totally real $2$-plane is $\mathrm{PU}(2,1)$-equivalent
to $\mathbf  {H}^2_{\mathbb R}$, the orthogonal projection to any
totally real $2$-plane enjoys the same properties.

Let $\langle H_r\rangle$ be the group of all
horizontal translation $H_r(x,y,u,v)=(x+r,y,u,v-2ry)$, and
let $\langle V_t\rangle$ be the group of all vertical translations
$V_t(x, y, u, v)=(x, y, u, v+t)$.
Thus $\langle V_t\rangle$ and $\langle H_r\rangle$
are Lie subgroups of $\mathcal H$, each isomorphic to the group 
of additive reals.

One can check that the quotient map 
$\Pi\co Y\to Y/\langle H_r\rangle$
is the smooth trivial principal $\langle H_r\rangle$-bundle, and 
the restriction 
$\Pi{|_P}\co P\to Y/\langle H_r\rangle$ of $\Pi$ to the plane 
\[P=\{(x,y,u,v)\in Y\co x=0\}\] is a diffeomorphism, so that
$Y/\langle H_r\rangle$ is 
diffeomorphic to $\mathbb {R}^2\times [0,\infty)$.

Since vertical and horizontal translations commute,
$\langle V_t\rangle$ naturally acts on $Y/\langle H_r\rangle$
so that the diffeomorphism $\Pi{|_P}$ is $\langle V_t\rangle$-equivariant.
Thus since $\langle V_t\rangle$ acts on $P$ properly, 
it acts properly on $Y/\langle H_r\rangle$.
In particular, if $U$ and $U^\prime$ are subsets 
of $Y$ such that 
$\Pi(U)$ and $\Pi(U^\prime)$ are precompact in $Y/\langle H_r\rangle$, 
then $\{t\in \mathbb R\co V_t(U)\cap U^\prime\neq\emptyset\}$
is precompact.  

For $v\in\mathbb R$, let $\Sigma_v$ be the totally real 
$2$-plane in $\mathbf {H}^2_{\mathbb C}$ given by
$\{(x, 0, u, v)\in \mathbf {H}^2_{\mathbb C}\}$.
Let $\pi_v\co \mathbf {H}^2_{\mathbb C}\to \Sigma_v $ 
be the orthogonal projection, and let $\bar\pi_v:
\mathbf {H}^2_{\mathbb C}\cup\partial_\infty\mathbf {H}^2_{\mathbb C}\to 
\Sigma_v\cup\partial_\infty\Sigma_v$ be the
extension of $\pi_v$ as above.
Since $\Sigma_v$ is $\langle H_r\rangle$-invariant, $\bar\pi_v$ is
$\langle H_r\rangle$-equivariant. One can check that for any 
horoball $B\subset \Sigma_v$ centered at $\infty$,
the set $\Pi(\bar\pi_v^{-1}(\Sigma_v\setminus B))$ is precompact
in $Y/\langle H_r\rangle$.
Finally, it is easy to see that
if $l\subset \Sigma_v$ is a geodesic passing through $\infty $
and $L=\pi_v^{-1}(l)$, then the 
restriction $\Pi|{_L}\co L\to \mathbf {H}^2_{\mathbb C}/\langle H_r\rangle$ 
of $\Pi$ to $L$ is a diffeomorphism. Thus $L$ can be thought of as
a smooth section of the bundle $\Pi$.

\section{Maskit combination theorem}
One of the common methods of producing new
discrete groups is the so-called combination
theorems. Here we state a combination theorem  
for groups acting by homeomorphisms
on an arbitrary topological space 
due to B.~Maskit~\cite[VII.A]{Mas}.

Let $X$ be a topological space and $\Gamma$
be a subgroup of $\mathrm {Homeo}(X)$.
A subspace $Y\subset X$ is called {\it precisely invariant}
with respect to a subgroup $H\le\Gamma$ if $Y$ is
$H$-invariant, and $\gamma(Y)\cap Y=\emptyset$
for $\gamma\in\Gamma\setminus H$.
A subset $F\subset X$ is called a {\it fundamental set} for the
$\Gamma$-action on $X$ if $F$ contains no $\Gamma$-equivalent 
points and intersects every $\Gamma$-orbit.
We say that $\Gamma$ acts {\it discontinuously} on $X$
if it has a fundamental set with nonempty interior.

Let $\Gamma_1$, $\Gamma_2$
be subgroups of $\mathrm {Homeo}(X)$ and let $J$ be a subgroup 
of $\Gamma_1\cap\Gamma_2$.  
Following Maskit, a pair $(X_1, X_2)$ of 
disjoint, nonempty, $J$-invariant subsets 
of $X$ is called {\it proper interactive} if,
for each $m\in \{1,2\}$ every element of
$\Gamma _m\setminus J$ maps $X_m$ into $X_{3-m}$, and 
for some $m\in \{1,2\}$
there is a point in $X_{3-m}$ that is not 
$\Gamma_m$-equivalent to any point of $X_m$. 

\begin{thm}{\rm (Maskit)}\it\ \label{mask}
Let $\Gamma_1$ and $\Gamma_2$ be
subgroups of $\mathrm {Homeo}(X)$, and let $J$ be a subgroup 
of $\Gamma_1\cap\Gamma_2$ such that  
there exists a proper interactive pair $(X_1, X_2)$.
Assume that for each $m\in \{1,2\}$
there is a fundamental set $\Phi_m$ for the $\Gamma_m$-action on $X$ 
such that $\Phi_m$ has nonempty interior, $\Phi_m\cap X_m$ is a 
fundamental set for the $J$-action on $X_m$,
and every element of $\Gamma_m$ maps $\Phi_m\cap X_{3-m}$ into $X_{3-m}$.\\
Then the group $\Gamma=\langle\Gamma_1, \Gamma_2\rangle$
is isomorphic to $\Gamma_1*_J\Gamma_2$,
and the set $\Phi=(\Phi_1\cap X_2)\cup (\Phi_2\cap X_1)$ 
is precisely invariant under the identity in $\Gamma$.
In particular, if $\Phi$ has nonempty interior, then $\Gamma $ 
acts discontinuously on $X$. 
\end{thm}

\section{Main construction}

For the rest of the paper we fix a closed orientable surface $M$
and a nontrivial element $\gamma\in\pi_1(M)$ represented by a simple closed
loop $\tilde\gamma$ that separates $M$ into two connected
compact surfaces $M_1$, $M_2$ with boundary $\tilde\gamma$. 
Thus $\pi_1(M)$ can be written as
the amalgamated product
$\pi_1(M_1)\ast_{\langle\gamma\rangle}\pi_1(M_2)$.

For $m\in \{1, 2\}$, we
identify the interior of $M_m$ with the 
finite volume hyperbolic surface  
$\mathbf {H}^2_{\mathbb R}/\Gamma_m$ where $\Gamma_m$ is
a noncocompact lattice in $\mathrm {SO}_0(2,1)$.
We can assume that $\Gamma_m$ has a fundamental polyhedron 
with finitely many sides, and exactly two sides passing through 
the point infinity in the upper half plane $\mathbf {H}^2_{\mathbb R}$.
We call these two sides {\it unbounded} and all the other sides
{\it bounded}. 
After removing some boundary points of the polyhedron,
we get a fundamental set $F_m$ for $\Gamma_m$ 
in $\mathbf {H}^2_{\mathbb R}$. 
Slightly abusing notations, we usually do not distinguish between 
$F_m$ and the original polyhedron, 
in particular, we often talk about sides of $F_m$.
 
Let $\Sigma_{v_1}$, $\Sigma_{v_2}$ be totally real $2$-planes as in
section~\ref{sec: vh}. 
The stabilizer of $\Sigma_{v_m}$ in $\mathbf {PU}(2,1)$
is isomorphic to $\mathbf {SO}(2,1)$ so we
can make $\Gamma_m$ act on $\mathbf  {H}^2_{\mathbb C}$ stabilizing
$\Sigma_{v_m}.$ 
Furthermore, we identify $F_m$ with a polyhedron in $\Sigma_{v_m}$ 
such that unbounded sides of $F_m$ pass through $\infty$ in the 
Siegel domain $\mathbf  {H}^2_{\mathbb C}$.

As the projection $\pi_{v_m}\co 
\mathbf  {H}^2_{\mathbb C}\to \Sigma_{v_m}$ is $\Gamma_m$-equivariant,
$\Gamma_m$ acts properly discontinuously on
$\mathbf  {H}^2_{\mathbb C}$ with a fundamental set $\Psi_m=\pi^{-1}_{v_m}(F_m)$.
Since $F_m$ is the intersection of finitely many halfplanes,
$\Psi_m$ is the intersection of finitely many {\it halfspaces}.
Each halfspaces is bounded by a smooth hypersurface which we call
{\it hyperplane}. We say a hyperplane is {\it unbounded} if
it passes through $\infty$. Otherwise, a hyperplane is {\it bounded}.

Let $d_m\in\Gamma_m$ be a parabolic element pairing the
unbounded sides of $F_m$. 
Then $d_m$ acts as a horizontal translation preserving $\Sigma_{v_m}$,
that is, $d_m=H_{r_m}$ for some $r_m\in \mathbb R$. 
We now assume that $r_1=r_2$, i.e. $d_1=d_2$, and we denote
this horizontal translation by $d$.  

Let $\beta_m\subset b_m\subset B_m$ be three concentric horoballs in 
the plane $\Sigma_{v_m}$  centered at $\infty$.  
For each $m$, choose $B_m$ small enough so that 
$\Sigma_{v_m}\setminus B_m$ contains all bounded sides of $F_m$ 
(and,  hence, $\pi_{v_m}^{-1}(\Sigma_{v_m}\setminus B_m)$
contains all bounded hyperplanes of $\Psi_m$).

By section~\ref{sec: vh}, 
$\Pi$ maps the complements of $\pi_{v_m}^{-1}(\beta_m)$ and
$\pi_{v_{3-m}}^{-1}(\beta_{3-m})$ onto precompact subsets 
of $Y/\langle H_r\rangle$, hence
we can choose $t=v_2-v_1$ so large that the subsets are disjoint. 
In fact, we can assume that the subsets
lie on the different sides of a properly embedded
hypersurface $H\subset Y/\langle H_r\rangle$ 
which becomes a linear half plane under the identification
$Y/\langle H_r\rangle\simeq \mathbb R^2\times [0,\infty)$.

The hypersurface $\Pi^{-1}(H)$ divides 
$\mathbf  {H}^2_{\mathbb C}$ into
two connected components. We let $X_1$ be the component
containing the bounded hyperplanes of $\Phi_2$
and let $X_2$ be the closure of the other component.
Thus $X_2$ contains all bounded hyperplanes of $\Phi_1$
and $X_1\cup X_2=\mathbf  {H}^2_{\mathbb C}$. 

Let $D_1=\Pi(\pi^{-1}_{v_2}(\beta_2))$ and 
$D_2=\Pi(\pi^{-1}_{v_2}(\Sigma_{v_2}\setminus b_2))$;
these are disjoint domains in 
$\mathbf  {H}^2_{\mathbb C}/\langle H_r\rangle$.
Consider the unbounded hyperplanes $S_m$, $d(S_m)$ of $\Psi_m$.
By section~\ref{sec: vh},
$S_m$ can be thought of as
a smooth section of the bundle 
$\Pi\co\mathbf  {H}^2_{\mathbb C}\to
\mathbf  {H}^2_{\mathbb C}/\langle H_r\rangle$.
Using bump functions, we construct a smooth section $S$
of the bundle whose restriction to $D_m$ is $S_m$.
Notice that $S$ splits $\mathbf  {H}^2_{\mathbb C}$ 
into two connected components;
we call the component that contains $d(S)$ 
a {\it halfspace associated to} $S$. Similarly, 
$d(S)$ splits $\mathbf  {H}^2_{\mathbb C}$ in two components and we call
the component that contains $S$ 
a {\it halfspace associated to} $d(S)$. 
Let $\Phi_m$ be the intersection of all the bounded halfspaces of
$\Psi_m$ and the halfspaces associated to $S$ and $d(S)$. 

\section{Discreteness}

In this section, we show that $\Phi_m$ is a fundamental set 
for the $\Gamma_m$-action on $\mathbf  {H}^2_{\mathbb C}$, 
and $\Phi=\Phi_1\cap\Phi_2$ is the fundamental set
for the $\Gamma$-action on $\mathbf  {H}^2_{\mathbb C}$, where
$\Gamma=\langle\Gamma_1,\Gamma_2\rangle$.

\begin{lem} For each $m\in \{1,2\}$, 
$\Phi_m$ is a fundamental set for $\Gamma_m$.
\end{lem}
\begin{proof}
The intersection of the halfspaces associated to $S$ and $d(S)$
is a fundamental set for the group generated by $d$. 
By construction,
$\Psi_m\subset \bigcup_{n\in \mathbb Z}d^n(\Phi_m)$.
Hence 
$\mathbf  {H}^2_{\mathbb C}=
\bigcup_{\gamma\in\Gamma_m}\Psi_m\subset
\bigcup_{\gamma\in\Gamma_m}\Phi_m$.
If $\Phi_m$ has $\Gamma_m$-equivalent points, they
must be $\langle d\rangle$-equivalent since
$\Phi_m \subset \bigcup_{n\in \mathbb Z}d^n(\Psi_m)$.
By construction, $\Phi_m$ has no $\langle d\rangle$-equivalent points
which completes the proof.
\end{proof}

\begin{lem} $(X_1, X_2)$ is a proper interactive pair.
\end{lem}
\begin{proof}
Clearly, $X_m$ is  $\langle d\rangle$-invariant, and 
$X_m\subset \bigcup_{n\in \mathbb Z}d^n(\Phi_m)$.
So for any $g\in \Gamma_m\setminus \langle d\rangle$,
\[g(X_m)\cap X_m\subset g(\bigcup_{n\in \mathbb Z}d^n(\Phi_m))
\cap \bigcup_{n\in \mathbb Z}d^n(\Phi_m)=\emptyset \]
since $\Phi_m$ is a fundamental set. 
Hence, $X_m$ is precisely invariant under
$\langle d\rangle$ in $\Gamma_m$.

Now take $g\in \Gamma_m\setminus \langle d\rangle$. 
Since $X_1$ is the complement of $X_2$, 
and $X_1$, $X_2$ are disjoint,
$g(X_m)\cap X_m=\emptyset$ implies  
$g(X_m)\subseteq X_{3-m}$ as desired.
It remains to find a point in $X_{3-m}$
that is not $\Gamma_m$-equivalent to any point of $X_m$.
Take $x\in \Phi_1\cap\Phi_2\cap X_{3-m}$
and assume $g(x)\in X_m$ for some $g\in \Gamma_m$.
Since $X_m\subset \bigcup_{n\in \mathbb Z}d^n(\Phi_m)$,
we get $d^ng(x)\in \Phi_m$ for some $n$.
Hence $g=d^{-n}$ and, therefore, $g(x)\in X_{3-m}$
because $X_{3-m}$ is $d$-invariant. Since $X_1, X_2$
are disjoint, we get a contradiction.
\end{proof}

\begin{thm} The group
$\Gamma=\langle\Gamma_1, \Gamma_2\rangle$ is discrete and isomorphic to
$\Gamma_1*_{\langle d\rangle}\Gamma_2\cong\pi_1(M)$. Moreover, 
$\Phi=\Phi_1\cap \Phi_2$ is a fundamental set
for $\Gamma$.
\end{thm}
\begin{proof} 
We want to apply Theorem~\ref{mask}.
First, obviously, $\Phi_m\cap X_m$ is a fundamental set for
the $\langle d\rangle$-action on $X_m$.
Second, we check that every element of $\Gamma_m$ maps
$\Phi_m\cap X_{3-m}$ into $X_{3-m}$.
(If not, then there is
$x\in g(\Phi_m\cap X_{3-m})\cap X_m$. Then for some $n$, 
$d^n(x)\in d^ng(\Phi_m\cap X_{3-m})\cap (\Phi_m\cap X_m)$.
Since $\Phi_m$ is a fundamental set, $g=d^{-n}$
which leads to a contradiction because $X_{3-m}$ is $d$-invariant).
Thus, since $\Phi_1\cap \Phi_2=(\Phi_2\cap X_1)\cup (\Phi_1\cap X_2)$ 
has nonempty interior
we can apply Theorem~\ref{mask}
to conclude $\Gamma$ is discrete, isomorphic to
$\Gamma_1*_{\langle d\rangle}\Gamma_2$ and 
$\Phi$ is precisely invariant under the identity in $\Gamma $.

We show that $\Phi=\Phi_1\cap \Phi_2$ is a fundamental 
set for $\Gamma $ by proving that the projection 
$p\co\mathbf  {H}^2_{\mathbb C}\to 
\mathbf  {H}^2_{\mathbb C}/\Gamma=N$ 
maps $\Phi$ onto an open and a closed subset of $N$, 
so that $p(\Phi)=N$.

We first show that $p(\Phi)$ is open. Since $p$ is an open map 
and $p(\Phi)=p(\bigcup_{g\in \Gamma}g(\Phi))$,
it is enough to prove that $\bigcup_{g\in \Gamma}g(\Phi)$
is open. Any point of the set is $\Gamma $-equivalent
to a point of $\Phi$, hence, it suffices to show that
$\Phi$ has an open neighborhood in 
$\bigcup_{g\in \Gamma}g(\Phi)$. 
One easily sees that
$\Phi$ lies in the interior of the set 
$\bigcup_{g\in S}g(\Phi)$ where 
$S=\{g\in \Gamma\co \overline{\Phi}\cap g(\overline{\Phi})
\neq\emptyset\}$.

Now prove that $p(\Phi)$ is closed in $N$.
Let $x_n\in \Phi$ be an arbitrary sequence such that
$p(x_n)$ converges to $y\in N$.
If $x_n$ subconverges to $z\in\overline{\Phi}$, then  
$p(z)=y$ and we are done since $p(\Phi)=p(\overline{\Phi})$.
It remains to consider the case when $x_n$ converges to $\infty$.
Passing to a subsequence, we can assume that
$x_n\in\pi_{v_2}^{-1}(\beta_2)\cap\Phi$, and  
there is a sequence $y_n$ with $p(y_n)=y$
and $\mathrm{dist}(x_n,y_n)<1$. 
Note that $\pi_{v_2}^{-1}(\beta_2)\cap\Phi=
\pi_{v_2}^{-1}(\beta_2)\cap\Psi_1\subset \pi_{v_1}^{-1}(F_1)$.
Hence $\pi_{v_1}(x_n)\in F_1$ and 
since $x_n$ converges to $\infty$, we can assume 
$x_n\in F_1\cap\beta_1$ and
$\mathrm{dist}(\pi_{v_1}(x_n),\partial \beta_1)>1$.
Since $\pi_{v_1}$ is distance
nonincreasing, $\pi_{v_1}(y_n)\in \beta_1$.
Thus, $y_n\in\pi_{v_2}^{-1}(\beta_1)\subset\cup_{n\in\mathbb Z} d^n(\Phi)$
so that $y\in p(\Phi)$.
\end{proof}

\section{Geometrical finiteness}
 
We refer~\cite{Bow} for background on geometrical finiteness
for manifolds of pinched negative curvature. 
\begin{thm} 
The group $\Gamma$ is a geometrically finite and any
maximal parabolic subgroup of $\Gamma$ is conjugate
to the cyclic subgroup generated by $\rho(\gamma)$. 
\end{thm}
\begin{proof} 
One of the definitions of 
geometrical finiteness given in~\cite{Bow}
is that the quotient manifold 
has finitely many ends, and every end is standard parabolic. 
Each of the manifolds $\Phi/\Gamma$ and $\Psi_1/\Gamma_1$
has exactly one end which is isometric to the quotient of
\[\pi^{-1}_{v_2}(\beta_2)\cap\Phi=\pi^{-1}_{v_2}(\beta_2)\cap\Psi_1\] 
by the subgroup $\langle d\rangle$.
The group $\Gamma_1$ is of course geometrically finite
as a subgroup of $\mathrm{SO}(2,1)$, hence 
$\Gamma_1$ is geometrically finite
as a subgroup of $\mathrm{PU}(2,1)$
since geometric finiteness is encoded in the $\Gamma$-action on
the limit set~\cite{Bow}. Therefore, $\Gamma $ is geometrically finite. 
\end{proof}

\begin{rmk}
Since $\Gamma $ is geometrically finite, the limit set 
$\Lambda(\Gamma )$ of $\Gamma $ consists of conical limit points and bounded
parabolic points. The set of bounded parabolic points is precisely
the $\Gamma$-orbit of $\infty$.
Being a closed surface group $\Gamma $ acts on 
$\mathrm{S}^1=\partial_{\infty}\mathbf  {H}^2_{\mathbb R}$. 
Now Floyd's theorem~\cite{Flo} (adapted for the complex hyperbolic case) implies
that there is a continuous $\Gamma $-equivariant map
of $\mathrm{S}^1$ onto $\Lambda(\Gamma )$ such that every conical limit point
has one preimage and every bounded parabolic point has exactly
two preimages.
\end{rmk}

\section{Quotients are tangent bundles}

\begin{thm} \label{thm: tb}
The quotient manifold $N=\mathbf  H^2_{\mathbb C}/\Gamma$
is diffeomorphic to the total space of the tangent bundle to $M$. 
\end{thm}

\begin{proof} 
Our first goal is to define a $\Gamma $-equivariant
$\mathbb R^2$-bundle structure on $\overline{\Phi}$ 
(in this proof all closures are taken in $\mathbf  H^2_{\mathbb C}$).
The set $\overline{\Phi}=\overline{\Phi_1}\cap\overline{\Phi_2}$
is the union of three pieces:
\[\overline{\Phi_1}\setminus \pi^{-1}_1(b_1)\!\ ,\ \
\overline{\Phi_2}\setminus \pi^{-1}_2(b_2)\ \ \mathrm{and}\ \
\pi^{-1}_1(B_1)\cap \pi^{-1}_2(B_2)\cap \overline{\Phi}\]
where 
$\overline{\Phi_1}\setminus \pi^{-1}_1(b_1)$ and
$\overline{\Phi_2}\setminus \pi^{-1}_2(b_2)$ are disjoint. 
Moreover, the intersection of $\overline{\Phi_m}\setminus \pi^{-1}_m(b_m)$
and $\pi^{-1}_1(B_1)\cap \pi^{-1}_2(B_2)\cap \overline{\Phi}$
is the set $\pi^{-1}_{v_m}(B_m\setminus b_m)\cap \overline{\Phi}$.
The map 
\[\pi_{v_m}\co\overline{\Phi_m}\setminus \pi^{-1}_m(b_m)\to
\overline{F_m}\setminus b_m\] 
is a smooth $\Gamma $-equivariant
$\mathbb R^2$-bundle whose fibers are totally real $2$-planes 
orthogonal to $\Sigma_{v_m}$. 
We now define a smooth $\Gamma $-equivariant
$\mathbb R^2$-bundle structure on 
$\pi^{-1}_1(B_1)\cap \pi^{-1}_2(B_2)\cap \overline{\Phi}$
in such a way that on the overlap 
$\pi^{-1}_{v_m}(B_m\setminus b_m)\cap \overline{\Phi}$
it coincides with 
\[\pi_{v_m}\co\pi^{-1}_{v_m}(B_m\setminus b_m)\cap \overline{\Phi}\to
(B_m\setminus b_m)\cap \overline{F_m}.\]
Find a smooth proper embedding 
$f\co\mathbb R\to \pi^{-1}_1(B_1)\cap \pi^{-1}_2(B_2)\cap S$
such that the intersection of $f(\mathbb R)$ and 
$\pi^{-1}_{v_m}(B_m\setminus b_m)\cap S$
is the interval $B_m\setminus b_m\cap S$.
It is easy to construct a smooth $\mathbb R^2$-bundle
$\pi^{-1}_1(B_1)\cap \pi^{-1}_2(B_2)\cap S\to f(\mathbb R)$ that
extends the bundle 
\[\pi_{v_m}\co\pi^{-1}_{v_m}(B_m\setminus b_m)\cap S\to 
(B_m\setminus b_m)\cap S.\]
This defines an $\langle H_r\rangle$-equivariant smooth $\mathbb R^2$-bundle
structure on $\pi^{-1}_1(B_1)\cap \pi^{-1}_2(B_2)$ over the
$\langle H_r\rangle$-orbit of $f(\mathbb R)$
that extends the bundle
$\pi_{v_m}\co\pi^{-1}_{v_m}(B_m\setminus b_m)\to B_m\setminus b_m$.

Let $E$ be the intersection of $\overline{\Phi}$
and the $\langle H_r\rangle$-orbit of $f(\mathbb R)$,
and let
$D=(\overline{F_1}\setminus b_1)\cup E\cup (\overline{F_2}\setminus b_2)$.
We have just constructed a smooth $\Gamma $-equivariant
$\mathbb R^2$-bundle $\overline{\Phi}\to D$. 
Let $\tilde D$ be the $\Gamma $-orbit of $D$. 

It is immediate to check that $\tilde D$ is a smooth submanifold
of $\mathbf  H^2_{\mathbb C}$ and the $\Gamma$-action on $\tilde D$
is smooth, free, and properly discontinuous so that $\tilde D/\Gamma$
is diffeomorphic to $M$. Therefore, we get
a smooth $\Gamma $-equivariant $\mathbb R^2$-bundle 
$\mathbf  H^2_{\mathbb C}\to \tilde D$. Passing to quotients
yields a smooth $\mathbb R^2$-bundle 
$N=\mathbf  H^2_{\mathbb C}/\Gamma\to \tilde D/\Gamma=M$.  

Choose orientations on $\Sigma_{v_1}$, 
$\Sigma_{v_2}$ so that $\Gamma_m$ preserves the orientation
on $\Sigma_{v_m}$ and the vertical translation
$V_t\co\Sigma_{v_1}\to \Sigma_{v_1}$ is orientation
preserving. Since $\Sigma_{v_m}$ is totally real,
the complex structure $\mathbb J$ defines a $\mathrm{SO}(2,1)$-equivariant
isomorphism between the normal and tangent bundles to 
$\Sigma_{v_m}$. The above orientation on 
$\Sigma_{v_m}$ together with its $\mathbb J$-image
defines an orientation on $\mathbf  {H}^2_{\mathbb C}$
which coincides with the canonical orientation of 
$\mathbf  {H}^2_{\mathbb C}$. 
This defines the orientations on $M$, $N$, and 
the bundle $N\to M$.

Oriented plane bundles over an oriented closed surface 
$M$ are (smoothly) isomorphic iff their Euler 
numbers are equal, so
it suffices to show that the Euler number of the bundle 
$N\to M$ is $\chi (M)$.
The Euler number is equal to the   
self-intersection number of the ``zero section'' 
$\sigma_0\co M=D/\Gamma\to N$. 
The surface $M=D/\Gamma$ is the union of three pieces
\[M_1=(\overline{F_1}\setminus b_1)/\Gamma\!\ ,\ \
T=E/\Gamma\ \ \mathrm{and}\ \
M_2=(\overline{F_2}\setminus b_2)/\Gamma.\]
Choose a smooth section $\sigma\co M\to N$ which is transverse to
$\sigma_0$.
Since $T$ is an annulus, we can assume that $\sigma (T)$ is disjoint from
$\sigma_0(M)$. So the intersection number of $\sigma $ and $\sigma_0$
is the sum of the intersection numbers $i_1$ and $i_2$ where
$i_m$ is intersection number of $\sigma $ and $\sigma_0$
restricted to $M_m$. By construction the bundle $N\to M$ restricted to $M_m$
is isomorphic to the tangent bundle of $M_m$ (as the universal
cover of $M_m$ is "totally real"). It is now a standard computation
to see that $i_m=\chi (M_m)$. (Hint: double $M_m$ along the boundary
to produce a closed oriented surface $S$ with $\chi (S)=2\chi (M_m)$. 
The intersection number of the double of ${\sigma_0}|_{M_m}$ and the double of
$\sigma |_{M_m}$ is $2i_m$. 
On the other hand, the self-intersection number of the zero section of 
$\mathrm{T}S$ is $\chi (S)$.)
Thus the self-intersection number of $\sigma_0$ is 
$i_1+i_2=\chi (M_1)+\chi (M_2)=\chi (M)$ as desired. 
\end{proof}

\section{Computing the Toledo invariant}

For background on Toledo invariant see~\cite{Tol, GKL, Xia}.
Here is a short account sufficient for our purposes.
Let $\rho\co\pi_1(M)\to \mathbf  {PU}(2,1)$ be a discrete faithful
representation. The K{\"a}hler form $\omega$ on $\mathbf  {H}^2_{\mathbb C}$
defines the K{\"a}hler form $\omega_N$ on $N=\mathbf  {H}^2_{\mathbb C}/\Gamma$.
Consider a smooth homotopy equivalence $f\co M\to N$ uniquely defined by $\rho$
up to homotopy. The Toledo invariant of $\rho$ is 
$\tau(\rho)=\frac{1}{2\pi}\int_M{f^*\omega_N}$. 
It is proved in~\cite{GKL}
that $\tau(\rho)\in\frac{2}{3}\mathbb Z$.

\begin{thm} 
If $\rho(\pi_1(M))=\Gamma$, then 
$\tau(\rho)=0$. 
\end{thm}

\begin{proof} 
The surface $M\setminus T$ in $N$ is
the union of two disjoint surfaces $M_1$ and $M_2$. 
The universal cover of each
of them lies in a totally real subspace. Since the K\"ahler form $\omega$
vanishes when restricted to a totally real subspace, we get
$\int_{M_m}{f^*\omega_N}=0.$
Therefore,
\[\tau(\rho)=
\frac{1}{2\pi}\int_{M_1}{f^*\omega_N}+
\frac{1}{2\pi}\int_{M_2}{f^*\omega_N}+
\frac{1}{2\pi}\int_{T}{f^*\omega_N}=
\frac{1}{2\pi}\int_{T}{f^*\omega_N}=
\frac{1}{2\pi}\int_{E}{\omega}.\]
Consider the horizontal translation $d=H_r$ that identifies
the boundary components of $E$. Consider the positive integer 
$n=3|\tau(\rho)|+1$, 
and let $s=\frac{r}{n}$.
Let $e$  and $H_r(e)$ be boundary curves of $E$.
Call $E_i$ the subsurface of $E$ bounded by $H_{is}(e)$ and
$H_s(H_{is}(e))$ where $i=0,\dots, n-1$. 
Since the form $\omega $ is 
$\mathbf {PU}(2,1)$-invariant, we get
\[\int\limits_{E_i}{\omega}=\!\!\!\!\!
\int\limits_{H_{is}(E_0)}{\!\!\!\!\omega}=
\int\limits_{E_0}H_{is}^*{\omega}=
\int\limits_{E_0}{\omega}.\]
Therefore,
\[\tau(\rho)=\frac{n}{2\pi}\int_{E_0}{\omega}\ \ \ \ 
\mathrm{so}\ \mathrm{that}\ \ \ \
\left|\frac{1}{2\pi}\int_{E_0}{\omega}\right|<\frac{1}{3}.\]
To prove $\tau(\rho)=0$, it suffices to show that 
$\frac{1}{2\pi}\int_{E_0}{\omega}\in\frac{2}{3}\mathbb Z$.
Notice that in the construction of 
$\Gamma =\Gamma_1*_{\langle d\rangle}\Gamma_2$ 
the parabolic element $d=H_r$ can be chosen arbitrarily. In particular,
one can take $r=s$. This defines a new 
discrete faithful representation  
$\rho_s\co \pi_1(M)\to\mathrm{PU}(2,1)$. 
Repeating the construction of Theorem~\ref{thm: tb}, 
one can choose the surface $E$ so that it coincides with the surface 
that was denoted $E_0$ above. 
Then $\frac{1}{2\pi}\int_{E_0}{\omega}$ is equal to the Toledo invariant
of $\rho_s$ which lies in $\frac{2}{3}\mathbb Z$ as needed.
\end{proof}

\small
\bibliographystyle{amsalpha}
\bibliography{ap}

\

DEPARTMENT OF MATHEMATICS, 253-37, CALIFORNIA INSTITUTE OF TECHNOLOGY,

PASADENA, CA 91125, USA

{\normalsize
{\it email:} \texttt{ibeleg@its.caltech.edu}}

\end{document}